\documentstyle[12pt,twoside]{article}
\newtheorem{theorem}{Theorem}[section]

\newtheorem{definition}{Definition}[section]
\newtheorem{proposition}{Proposition}[section]
\newtheorem{corollary}{Corollary}[section]
\newtheorem{lemma}{Lemma}[section]

\newcommand{\proof}{\noindent{\it Proof. }}
\newcommand{\bbfR}{{I\!\!R}}
\newcommand{\bbfN}{{I\!\!N}}
\renewcommand{\P}{{I\!\!P}}

\renewcommand{\div}{\nabla\cdot}
\renewcommand{\S}{{\cal S}}
\newcommand{\C}{{\cal C}}

\newcommand{\rf}[1]{(\ref{#1})}

\newcommand{\X}{{\cal X}}

\def\cbdu{\hfill{$\Box$}}

\def\ya{{\cal Y}^a}
\def\x{{\cal X}}
\def\esssup{\mbox{\rm ess} \sup_{\!\!\!\!\!\!\!\!\! \xi\in\bbfR^3}}
\newcommand{\PM}{{{\cal P\! M}}}
\newcommand{\DPM}{{\PM^2}}

\newcommand{\PMa}{{ {\cal P\!M}^a}}
\newcommand{\I}{{|\!|\!|}}

\newcounter{remark}
\newenvironment{remark}%
{\medskip \stepcounter{remark} \noindent {\it Remark
\arabic{section}.\arabic{remark}} }{\rm \cbdu}

\title{\bf
Smooth  or singular solutions to the Navier--Stokes system ?}

\author{{\sc Marco Cannone \& Grzegorz Karch}\\
\\
{\small Universit\'e de Marne-la-Vall\'ee, Laboratoire d'Analyse}\\
{\small et de Math\'ematiques Appliqu\'ees}\\
{\small Cit\'e Descartes--5, bd Descartes, Champs-sur-Marne,}\\
{\small 77454 Marne-la-Vall\'ee Cedex 2, France}\\
{\small {\tt cannone@math.univ-mlv.fr}}\\
    \\
{\small Instytut Matematyczny, Uniwersytet Wroc\l awski}\\
{\small pl. Grunwaldzki 2/4, 50-384 Wroc\l aw, Poland;}\\
{\small Mathematical Institute}\\
{\small Polish Academy of Sciences, Warsaw
(2002-2003)}\\
{\small {\tt karch@math.uni.wroc.pl}}\\
}

\begin{document}
\maketitle

\noindent

{\parskip = 4pt plus 4pt minus 1pt}

\begin{abstract}
The existence of singular solutions of the incompressible Navier-Stokes 
system with  singular external forces, the
existence of
regular solutions for more regular forces as well as the asymptotic
stability of small solutions (including stationary ones), and a
 pointwise loss of smoothness for
solutions  are proved in the same function
space
of pseudomeasure type.


\end{abstract}

\newpage

\baselineskip=16pt

\section{Introduction}

So far, only two ways for attacking the Cauchy problem
for the Navier--Stokes equations are
known:
the first is due to J. Leray \cite {Ler}, and the second is due to T.
Kato \cite{Ka}.
None of them can be considered the
``golden rule'' for solving the Navier--Stokes equations because  they
both
leave open the following celebrated question.  In three dimensions,
does the
velocity field of a fluid flow that starts smooth remain smooth and
unique for
all time ?

The concept of ``weak'' solutions introduced by J. Leray in 1933,
permits the study of functions in
much larger classes than the classical spaces  used to describe
the motion of a fluid. It is easier to prove the
existence of a solution  (regular or singular) in a larger class, but
such a solution may
not be unique.  Based on {\sl a priori} energy estimates, Leray's
theory
gives the existence of global weak,
possibly  irregular and possibly non-unique solutions to the
Navier--Stokes equations.
On the other hand, a completely different theory introduced by T.
Kato in 1984, based on
semigroups techniques and the
fixed point scheme,
gives the existence of a global unique regular ``mild'' solution,
under the restrictive assumption of small
initial data. A second restriction is given by the fact that Kato's
algorithm  does not provide a framework for
studying {\sl a~priori} singular solutions. In fact, in order to
overcome
the difficulty (and sometimes the impossibility) of proving the
continuity of
the bilinear estimate in the, so-called, critical spaces, Kato's
algorithm makes clever use of a combination
of two estimates in two different norms, the natural one and a
regularizing norm. As such,
Kato's approach imposes
{\sl a~priori} a~regularization effect on solutions we look for. In
other
words, they are considered as fluctuations around the solution of the
heat equation with same initial data. For people who believe in blow
up and singularities,
this {\sl a~priori} condition coming from the ``two norms approach''
is indeed very strong.
However, there exist two exceptions, more exactly two critical spaces
where Kato's method applies with
just  one norm: the Lorentz space $L^{3,\infty}$ (considered
independently by  M. Yamazaki \cite{Y00}
and  by Y. Meyer \cite{M}) and the pseudomeasure space of Y. Le Jan
and A.S. Sznitman \cite{LJS}, \cite{CP2}.
Here we will not go into the technical details arising from these
critical spaces and we refer the reader
to the recent surveys contained in \cite{CaH} and in \cite{LR}.

    In this paper we will show how the approach with only one norm
gives existence and uniqueness of a
(small) solution in a larger space which, in our case, contains
genuinely singular solutions
that are not smoothed out by the action of the nonlinear
semigroup associated. More exactly, in the case of the pseudomeasure
space we can
prove the following results. The existence of singular solutions
associated to
singular (e.g. the Dirac delta) external forces thus  allowing to
describe the solutions
considered  by L.D. Landau in \cite{L} and by G. Tian and Z. Xin in
\cite{TX}. The
existence of regular solutions for more regular
external forces.  The asymptotic stability of small solutions
including
stationary ones. A
pointwise loss of smoothness for
solutions.

The study of the Navier--Stokes equations written in terms of the
vorticity and with measures
as initial data started in the 80s in a series of
papers  by G. Benfatto, R. Esposito, M. Pulvirenti \cite{BEP}, G.-H.
Cottet and J. Soler
\cite{Co, CoS}, and Y. Giga, T. Miyakawa and H. Osada \cite{GM, GMO}.
We refer
the reader as well to the more recent results obtained by T. Kato in
\cite{KaV} and Y. Giga in
\cite {GiV}. On  the other hand, the case of external forces that can
be singular atomic measures was studied by H. Kozono and M. Yamazaki
\cite{KY95}.
Here we want to provide, among others, such kind of results.


\section{One-point singular solutions}
\setcounter{equation}{0}
\setcounter{remark}{0}

As observed by J. Heywood in \cite{Hey}, in principle ``{\sl  it is
easy to construct a singular
solution of the NS equations that is driven
by a singular force. One simply constructs
a solenoidal vector field $u$ that begins
smoothly and evolves to develop a singularity,
and then defines the force to be the residual.}'' In this section we
want to give an explicit
example of this mathematical evidence. Our  example arises from
the physical experiment described by L.D. Landau in \cite{L} (see
also \cite[Sec. 23]{LL}), where an axially symmetric jet
discharging from a thin pipe into the  unbounded space is
studied. Passing to the limit with the diameter of the pipe,
this ``plunged'' jet can be regarded as emerging from a point source (i.e.
driven by the delta function). Landau provided a mathematical  setting for
explaining this phenomenon by using the classical
incompressible Navier--Stokes system  and deriving an explicit
``solution'' for it.

To be more precise, let us recall the famous Navier--Stokes equations,
describing the evolution
of the velocity field $u$ and pressure field $p$ of a
three-dimensional incompressible viscous fluid
at time $t$ and the position
$x\in\bbfR^3$. These equations are given by
\begin{eqnarray}
&&u_t-\Delta u +(u\cdot \nabla) u +\nabla p =F,\label{eq}\\
&&\div u =0, \label{div}\\
&&u(0)=u_0.\label{ini}
\end{eqnarray}
    where the external
force $F$ and initial velocity $u_0$ are assigned.

Recently, G. Tian and Z. Xin \cite{TX} also found explicit formulas
for a one-parameter family of stationary ``solutions''  of the
three-dimensional Navier--Stokes system ``with $F\equiv 0$'' which are
regular except at a given point. Due to the translation invariance of
the Navier--Stokes system, one can assume that the
singular point corresponds to the origin.
These explicit ``solutions'' by Tian and
Xin agree with those obtained by Landau for
special values of the parameter. More exactly, the main theorem from
\cite{TX} reads as follows. {\it All solutions
to system \rf{eq}--\rf{ini} (with $F\equiv 0$) $u(x)=(u_1(x), u_2(x),
u_3(x))$ and $p=p(x)$ which are steady,
symmetric about $x_1$-axis, homogeneous of degree $-1$, regular
except $(0,0,0)$ are given by the following explicit formulas:
\begin{eqnarray}
u_1(x)= 2{c|x|^2-2x_1|x|+cx_1^2\over |x|(c|x|-x_1)^2},  &&\quad
u_2(x)= 2{x_2(cx_1-|x|) \over |x|(c|x|-x_1)^2}, \label{sing-sol} \\
u_3(x)= 2{x_3(cx_1-|x|)  \over |x|(c|x|-x_1)^2},&&\quad
p(x)= 4{cx_1-|x|  \over |x|(c|x|-x_1)^2} \nonumber
\end{eqnarray}
where $|x|=\sqrt{x_1^2+x_2^2+x_3^2}$ and $c$ is an arbitrary constant
such that $|c|>1$.
}

\begin{remark}
Note that  in the formula \cite[(2.1)]{TX} the numerator of the
fraction defining $u_1(x)$ should read 
$cr^2-2r(x_1-x_1^0)+c(x_1-x_1^0)^2$. The factor ``$2$'' was missing
in that formula what can be inferred from \cite[(2.40)]{TX} or
\cite[(23,16)--(23,19)]{LL}. On the other hand, the sign ``$-$'' 
 in the formula \cite[(23,20)]{LL} for the pressure is wrong.
\end{remark}
\medskip

Before commenting this result, we think it is necessary to clarify
the meaning of ``solution of the
Navier--Stokes  equations'', for, since the appearance of the pioneer
papers of Leray, the word ``solution''
has been  used in a more or less generalized sense giving origin to
so many different definitions of
``solutions'', distinguished only  by the class of functions  they
are supposed to belong to: {\sl classical,
strong, mild,  weak, very weak, uniform weak} and  {\sl local Leray}
solutions of the Navier--Stokes equations ! We
will not present all the possible (more or less well-known)
definitions here and refer the reader to \cite{CaH}
and the references therein.

Let us first remark that there is no hope to describe the
``solutions''given by equations \rf{sing-sol} in
Leray's theory, because they are not globally of finite energy, in
other words
they do not belong to $L^2(\bbfR^3)$. However, they do belong to
$L^2_{loc}(\bbfR^3)$  and this is
at least enough to allow us to
give a (distributional) meaning to the nonlinear term
$(v\cdot\nabla)v=\div (v\otimes v)$.
Moreover, the ``solutions'' discover by Tian and Xin cannot be
analyzed by Kato's two norms
method either, because they are global but not smooth,  more exactly
they
are singular at the origin with a singularity of the kind $\sim
1/\vert x\vert$ for all time.

We will provide in the following section an {\sl ad hoc} framework
for
studying
such singularity within the fixed point scheme and without
using the two norms approach. As recalled in the introduction, this
can be done in principle either
in a Lorentz or in a pseudomeasure space and they both contain
singularities of the
type $\sim 1/\vert x\vert$. However, we will chose the latter space
not only because the
proofs will be very elementary, but also because this choice will
allow us to treat singular (Delta type)
external force, that precisely arise from Landau and Tian and Xin
``solutions''.

More exactly, by straightforward calculations, one can check that,
indeed, the functions
$(u_1(x), u_2(x), u_3(x))$ and $p(x)$ given by \rf{sing-sol} satisfy
\rf{eq}--\rf{ini} with
$F\equiv 0$ in the {\sl pointwise sense} for every $x\in
\bbfR^3\setminus
\{(0,0,0)\}$. On the other hand, if one treats $(u(x),p(x))$ as a
{\sl distributional or generalized} solution to
\rf{eq}--\rf{ini} in the whole
$\bbfR^3$, they
correspond to the very singular external force
$F=(b\delta_0, 0,0),$
where the parameter $b\neq 0$ depends on $c$ and $\delta_0$ stands for the
Dirac delta. Let us state this fact more precisely.

\begin{proposition}\label{prop:sing-sol}
Let $u=(u_1,u_2,u_3)$ and $p$ be defined by \rf{sing-sol}. For every
test function $\varphi \in C^\infty_c (\bbfR^3)$ the following
equalities hold true:
\begin{equation}
\int_{\bbfR^3} u\cdot \nabla \varphi \;dx =0 \label{w-div}
\end{equation}
and
\begin{equation}
\int_{\bbfR^3} \left(\nabla u_k\cdot \nabla \varphi -
u_ku\cdot \nabla \varphi -p {\partial\over \partial x_k} \varphi\right) \;
dx=
\left\{
\begin{array}{lcl}
-b(c)\varphi(0) &\mbox{if} & k=1\\
0&\mbox{if} & k=2,3,
\end{array}
\right.\label{w-eq}
\end{equation}
where
\begin{equation}b(c) =4\pi \left(4c+2c^2\log{c-1\over c+1} +{16\over
3}{c\over
c^2-1}\right).\label{fact.b}
\end{equation}
In particular, the function $b=b(c)$ is decreasing on $(-\infty, -1)$
and $(1,+\infty)$.
Moreover, $\lim_{c \searrow 1} b(c)=\infty$, $\lim_{c \nearrow -1}
b(c)=-\infty$ and
$\lim_{\vert c\vert \to\infty } b(c)=0$.
\end{proposition}

\proof
Equality \rf{w-div} says that the velocity $u$ is weakly
divergence-free in $\bbfR^3$. This can be shown by a standard
argument
involving integration by parts, since each component of $u$ is
homogeneous of degree $-1$ and thus belongs to $W^{1,p}_{loc} (\bbfR^3)$
with $1\leq p<3/2$
and $(\nabla \cdot u)(x)=0$ for all $x\in \bbfR^3\setminus \{0\}$.

Next, due to singularities of $u$ and $p$ at the origin, we fix
$\varepsilon>0$
and we integrate in equations \rf{w-eq} for $|x|\geq \varepsilon$,
only. Integrating by parts, we obtain
\begin{eqnarray}
&&\hspace{-1.5cm}
\int_{|x|\geq \varepsilon} \left(\nabla u_k\cdot \nabla \varphi -
u_ku\cdot \nabla \varphi -p {\partial\over \partial x_k} \varphi \right)\;
dx\nonumber\\
&=&\int_{|x|\geq \varepsilon}
\left(-\Delta u_k +\nabla\cdot (u_ku) +{\partial \over \partial x_k}p
\right)\varphi \;dx\label{x:eps}\\
&&+ \int_{|x|=\varepsilon} \left(
(\nabla u_k -u_ku )\cdot {x\over \varepsilon} -p {x_k\over
\varepsilon}\right) \varphi \;d\sigma(x),\nonumber
\end{eqnarray}
because $x/\varepsilon$ is the unit  vector normal to the sphere
$\{x\in\bbfR^3\,:\, |x|=\varepsilon\}$. Obviously, the first term on
the right-hand side of \rf{x:eps} disappears, and our goal is to
compute the limit as $\varepsilon \searrow 0$ of the second one.

For this reason, note first that each term $\nabla u_k$, $u_ku$, and
$p$ is
homogeneous of degree $-2$. Hence, changing variables $x=\varepsilon
y$ in the integral $\int_{|x|=\varepsilon} ...\;d\sigma(x)$ in
\rf{x:eps}, and next passing to the limit with $\varepsilon\searrow
0$ we show by the Lebesgue Dominated Convergence Theorem that it
converges toward
\begin{equation}
\varphi(0) \int_{|x|=1} \left( (\nabla u_k -u_ku)\cdot x -px_k\right)
\;d\sigma(x).
\label{x:eps2}
\end{equation}

To complete this proof, it remains to compute the surface integral
in \rf{x:eps2}. First, however, we simplify it a little by using the
Euler theorem for homogeneous functions which in this case gives
$x\cdot \nabla u_k=-u_k$. Moreover, it follows from the definition of
$u_k$ and $p$ that
$$
u_1={1\over 2} px_1+{2\over c|x|-x_1}, \quad
u_2={1\over 2} px_2,\quad
u_3={1\over 2} px_3.
$$
Consequently, for $k=2,3$, the integral in \rf{x:eps2} equals 
$$
-\int_{|x|=1} \left( u_k+u_k (u\cdot x) +2u_k\right)\; d\sigma(x) =0,
$$
because $u_2$ and $u_3$ are odd functions with respect to $x_2$ and
$x_3$, respectively, and $u\cdot x$ is even.
In case of $k=1$, we use the identities
$$
u_1(x) =c+(c^2-1) \left( {c\over (c-x_1)^2}-{2\over c-x_1}\right)
\quad\mbox{and} \quad px_1=2u_1-{4\over c-x_1}
$$
valid for $|x|=1$, and the polar coordinates to show that
\begin{eqnarray*}
&&\hspace{-0.5cm}\int_{|x|=1} \left(u_1+u_1 (u\cdot x) +2u_1 -{4\over
c-x_1}\right)\;d\sigma(x)\\
&&= 2 \pi \int_{-1}^1 2\left(\left(c+(c^2-1)\left({c\over (c-x_1)^2}-{2\over c-x_1}
\right)\right) \left( 1+2{c^2-1\over (c-x_1)^2}\right) -{2\over c-x_1}\right)\;dx_1\\ 
&&=b(c).
\end{eqnarray*}
Here, we skip these long but rather elementary calculations.
\cbdu

\begin{remark}
As we have already emphasized, the stationary solutions defined in
\rf{sing-sol} are singular with singularity of the kind
${\cal O}(1/|x|)$ as $|x|\to 0$. This is the critical singularity in the
context of Proposition \ref{prop:sing-sol}, because as it was shown
by H.J.Choe and H.Kim \cite{CK}, every pointwise stationary solution
to system
\rf{eq}--\rf{ini} with $F\equiv 0$ in
$B_R\setminus\{0\}=\{x\in\bbfR^3\,:\, 0<|x|<R\}$ satisfying
$u(x)=o(1/|x|)$ as $|x|\to 0$ is also a solution in the sense of
distributions in the whole $B_R$. Moreover, it is shown in
\cite{CK} that under the additional assumption $u\in L^q(B_R)$ for
some $q>3$, then the stationary solution $u(x)$ is smooth in the
whole ball $B_R$. In other words, if $u(x)=o(1/|x|)$ as $|x|\to 0$
and $u\in L^q(B_R)$ for some $q>3$, then the singularity at the origin is
removable.
\end{remark}
%
%

\section{Definitions and spaces}
\setcounter{equation}{0}
\setcounter{remark}{0}

We will study global-in-time
solutions $u=u(x,t)$
to the Cauchy problem in $\bbfR^3$  for the
incompressible Navier--Stokes equations \rf{eq}--\rf{div}. As far as
$u=u(x,t)$ is a sufficiently
regular function, the equations
\rf{eq}--\rf{div}
can be rewritten as
$$
u_t-\Delta u +\nabla \cdot (u\otimes u) +\nabla p =F, \quad \div u
=0.
$$
If we recall that  the Leray projector on solenoidal vector
fields is given
by the formula
\begin{equation}
\P v = v-\nabla \Delta^{-1} (\nabla \cdot v)\label{proj}
\end{equation}
for sufficiently smooth functions $v=(v_1(x), v_2(x), v_3(x))$,
we formally transform the system \rf{eq}--\rf{div}  into
$$
u_t -\Delta u +\P \nabla \cdot (u\otimes u)=\P F, \quad \nabla \cdot
u=0.
$$
Finally, let us emphasize that we shall study the problem
\rf{eq}--\rf{ini}
{\sl via} the following  integral equation obtained from the Duhamel
principle
\begin{eqnarray}
u(t)&=& S(t) u_0 - \int_0^t  S(t-\tau) \P \nabla \cdot (u\otimes
u)(\tau)\;d\tau \label{duhamel}\\
& & +\; \int_0^t S(t-\tau) \P F(\tau)\;d\tau, \nonumber
\end{eqnarray}
where $S(t)$ is the heat semigroup given as the convolution with the
Gauss--Weierstrass kernel: $G(x,t)=(4\pi t)^{-3/2} \exp(-|x|^2/(4t))$.
To give a meaning to the Leray projector $\P$
(defined in \rf{proj}), let us first
recall that the Riesz transforms $R_j$ are the pseudodifferential
operators
defined in the Fourier variables as $\widehat{R_k f} (\xi)
={i\xi_k\over |\xi|}
\widehat f(\xi)$. Here and in what follows the Fourier transform
of an integrable function $v$ is given by
$\widehat v(\xi)\equiv (2\pi)^{-n/2}\int_{\bbfR^n}
e^{-ix\cdot \xi} v(x)\;dx$. Using these well-known operators we
define
$$
(\P v)_j =v_j +\sum_{k=1}^3 R_j R_k v_k;
$$
moreover, in our considerations below, we shall often denote by
$\widehat \P (\xi)$
the symbol of the pseudodifferential operator $\P$ which is the
matrix
with components
$$
(\widehat \P (\xi))_{j,k} =\delta_{jk} -{\xi_j\xi_k\over |\xi|^2}.
$$
All these components are bounded on $\bbfR^3$ and we put
\begin{equation}
\kappa = \max_{1\leq j,k\leq 3} \sup_{\xi\in \bbfR^3 \setminus \{0\}}
|(\widehat \P (\xi))_{j,k}|.
\label{kappa}
\end{equation}

We are now in a position to introduce the Banach functional spaces
relevant to our study
of solutions of the Cauchy problem for the system \rf{eq}--\rf{ini}:
$$
\PM^a\equiv\{v\in{\cal S}'(\bbfR^d): \widehat v\in L^1_{\rm loc}(\bbfR^d),
\|v\|_{\PM^a}\equiv \esssup |\xi|^a|\widehat
v(\xi)|<\infty\},
$$
where $a\geq 0$ is a given parameter.
The notation $\PM$ stands for {\sl pseudomeasure}, and the
classical space of pseudomeasures
introduced in harmonic analysis (i.e.
those distributions whose Fourier transforms are bounded)  corresponds
to $a=0$.

\begin{definition}\label{def1}
By a solution of \rf{eq}--\rf{ini} we mean in this paper a function
$u=u(t)=(u_1(t), u_2(t), u_3(t))$ with each component $u_i$
belonging to the space of vector-valued functions
$\X={\cal C}_w([0,T);\PM^2)$, $0<T\le \infty$, and such that
\begin{eqnarray}
\widehat u(\xi,t)&=&e^{-t|\xi|^2}\widehat
u(\xi,0)+\int_0^te^{-(t-\tau)|\xi|^2}
\widehat \P (\xi)\; i\xi\cdot \left(\widehat{u\otimes
u}\right)(\xi,\tau)\, d\tau\label{FDuh}\\
&&+\int_0^te^{-(t-\tau)|\xi|^2}
\widehat \P (\xi) \widehat F(\xi, \tau)\, d\tau\nonumber
\end{eqnarray}
   for all $0\leq t\leq T$.
\end{definition}

The space $\PM^{2}$ is chosen because it contains  homogeneous
functions
of degree $-1$ which are sufficiently regular on the unit sphere.
In particular, one can easily check that this is the case for the
one-point singular solutions defined in \rf{sing-sol}.

\begin{remark}
Given $f\in \S'(\bbfR^3)\cap L^1_{loc}(\bbfR^3)$ we denote the
rescaling $f_\lambda (x)=f(\lambda x)$. In a~standard way, we extend
this definition to all tempered distributions.  It follows from
elementary calculations that $\widehat f_\lambda(\xi)=\lambda^{-3}
\widehat f(\lambda^{-1} \xi)$. Hence, for every $\lambda >0$, we
obtain the scaling property of the norm in $\PMa$
\begin{equation}
\| f(\lambda\,\cdot)\|_\PMa =\lambda^{a-3}\|f\|_\PMa.
\label{PMa-scal}
\end{equation}
In particular, the norm $\PM^2$ is invariant under rescaling
$f\mapsto\lambda f(\lambda\,\cdot)$.
Moreover, it follows from \rf{PMa-scal} that for $a=3(1-1/p)$ the
norms $\|\cdot\|_{\PMa}$ and $\|\cdot\|_{L^p(\bbfR^3)}$ have the same
scaling property.
\end{remark}

\begin{remark} ${\cal C}_w$ denotes, as usual (cf.
\cite{C}), the space of vector-valued functions which are weakly
continuous as distributions in $t$. This is an additional difficulty
caused by the fact that the heat semigroup $(S(t))_{t\ge 0}$
is not strongly continuous on the spaces of pseudomeasures but only
weakly continuous (cf. Lemma \ref{lem-a}, below).
\end{remark}

\begin{remark}
Usually, a {\sl mild solution} of an evolution equation like
\rf{eq}--\rf{ini}
is defined as a solution to the integral equation \rf{duhamel}
    and the integral is
understood as the Bochner integral. However, such a meaning of a
solution is not suitable for our construction of solutions of the
Cauchy problem and, in particular, of self-similar solutions.
Indeed, for stationary and homogeneous of degree $-1$ solutions $u$
(given, e.g., by \rf{sing-sol}), the nonlinear term  corresponds to
a tempered distribution which is homogeneous of degree $-3$, hence,
there exists a  distribution $H$ such that
$$
S(t) \P\nabla\cdot (u\otimes u) = t^{-3/2} H\left({\cdot\over \sqrt
t}\right).
$$
Now, computing the $\PM^2$ norm and using the scaling relation
\rf{PMa-scal}, we obtain
$$
\|S(t) \P\nabla\cdot (u\otimes u)\|_{\PM^2} =t^{-1} \|H\|_{\PM^2}.
$$
So, $S(t) (\P\nabla\cdot (u\otimes u))$ is not Bochner integrable as
a
mapping on $[0,T)$ with values in $\PM^2$.
On the other hand, the Fourier transform of this quantity equals to
$ e^{-t|\xi|^2} \widehat \P (\xi) (\widehat{u\otimes u})(\xi)$
and the singularity at $t=0$ does not appear. Hence, the integral
with respect to $\tau$  in equation \rf{duhamel} should be defined
in a weak sense like, e.g., it was done in \cite[Def. 2]{Y00}. 
For more explanations, 
we refer the reader to \cite{LR}, because our  spaces $\PM^a$
are the example of the shift-invariant Banach spaces of distributions
systematically used in that book.
\end{remark}

\medskip

Nevertheless, a distributional solution of system \rf{eq}--\rf{ini}
is a solution
of the integral equation of \rf{FDuh}, and vice versa. This
equivalence
can be proved by a standard reasoning, and we refer the interested
reader to \cite[Th. 5.2]{Y00} for details of such computations.
\bigskip

To simplify the notation, the quadratic term in
\rf{duhamel}
will be denoted  by
\begin{equation}
B(u,v)(t)=-\int_0^t S(t-\tau) \P \nabla \cdot (u\otimes
v)(\tau)\;d\tau,\label{Bf}
\end{equation}
where $u=u(t)$ and $v=v(t)$ are functions defined on $[0, T)$ with
values in a vector space (here most frequently $\PM^2$).


%
%

\section{Global-in-time solutions}
\setcounter{equation}{0}
\setcounter{remark}{0}

As in \cite{C}, the proof of our basic theorem
on the existence, uniqueness and stability of solutions to the
problem \rf{eq}--\rf{ini} is based on the following abstract lemma,
whose slightly more general form is taken
from \cite{LR}.

\begin{lemma}
\label{lem:xyB}
Let $(\X, \|\cdot\|_\X)$ be a Banach space and $B:\X\times
\X\to
\X$ a bounded bilinear form satisfying
$
\|B(x_1,x_2)\|_\X\leq \eta \|x_1\|_\X \|x_2\|_\X
$ for all $x_1,x_2\in\X$ and a constant $\eta>0$.
Then, if $0<\varepsilon<1/(4\eta)$ and if $y\in\X$  such that
$ \|y\|<\varepsilon$, the equation
$x=y+B(x,x)$ has a solution in $\X$ such that $\|x\|_\X\leq
2\varepsilon$.
This solution is the
only one in the ball $\bar B(0,2\varepsilon)$.
Moreover, the solution depends continuously on $y$ in the following
sense: if $\|\tilde y\|_\X\leq \varepsilon$,
$\tilde x=\tilde y+B(\tilde x, \tilde x)$, and $\|\tilde x\|_\X\leq
2\varepsilon$, then
$$
\|x-\tilde x\|_\X \leq {1\over 1-4\eta \varepsilon}\|y-\tilde y\|_\X.
$$
\end{lemma}

\proof Here, the reasoning is based on the standard Picard iteration
technique completed by the Banach fixed point theorem. For other
details of the proof, we refer the reader to \cite[Th. 13.2]{LR}.
\cbdu
\medskip

Our goal is to apply Lemma \ref{lem:xyB} in the space
\begin{equation}
\X={\cal C}_w([0,\infty), \PM^2)\label{X}
\end{equation}
to the integral equation \rf{duhamel}
which has the form $u=y+B(u,u)$, where the bilinear form is defined
in \rf{Bf}
and $y=S(t) u_0 +\int_0^t S(t-\tau) \P F(\tau)\;d\tau$.
We need some preliminary estimates.

\begin{lemma}\label{lem-a}
Given $u_0\in \PM^2$, we have $S(\cdot )u_0\in \X$.
\end{lemma}

\proof
By the definition of the norm in $\PM^2$, it follows that
$$
\|S(t) u_0\|_{\PM^2} =\esssup |\xi|^2\left|e^{-t|\xi|^2}
\widehat u_0(\xi)\right| \leq \esssup |\xi|^2|\widehat
u_0(\xi)|=\|u_0\|_{\PM^2},
$$
so, $S(\cdot)u_0\in L^\infty ([0, \infty), \PM^2)$.

Now, let us prove the weak continuity with respect to $t$, and, by
the  semigroup property of $S(t)$, it suffices to do this for $t=0$
only.
For every $\varphi\in \S(\bbfR^3)$, by the Plancherel formula, we
obtain
\begin{eqnarray*}
\left|\langle S(t)u_0-u_0, \varphi\rangle\right| &=&
\left|\int \left(e^{-t|\xi|^2}-1\right) \widehat u_0 (\xi) \widehat
\varphi(\xi)\;d\xi\right|\\
&\leq & t \; \esssup \left|{ e^{-t|\xi|^2}-1\over t|\xi|^2} \right|
\;\|u_0\|_{\PM^2}
\|\widehat\varphi\|_{L^1{\bbfR^3}} \to 0 \quad \mbox{as} \quad
t\searrow 0.
\end{eqnarray*}
\cbdu

\begin{lemma}\label{lem-b}
Given $F\in {\cal C}_w([0,\infty), \PM)$, it follows that
$$
w\equiv \int_0^t S(t-\tau) \P F(\tau) \;d\tau \in  \X.
$$
Moreover, $\|w\|_{\X}\leq \|F\|_{\C_w([0,\infty),
\PM)}$.
\end{lemma}

\proof
Similarly as in the proof of Lemma \ref{lem-a} we get
\begin{eqnarray*}
    \|w(t)\|_{\PM^2}&=& \esssup |\xi|^2 \left|
\int_0^t e^{-(t-\tau)|\xi|^2} \widehat P(\xi) \widehat F(\xi,
\tau)\;d\tau\right|\\
&\leq & \kappa |\xi|^2 \int_0^t e^{-(t-\tau)|\xi|^2}\;d\tau
\|F\|_{\C_w([0,\infty), \PM)}\\
&\leq & \kappa \|F\|_{\C_w([0,\infty), \PM)}.
\end{eqnarray*}
Let us skip the proof of the weak continuity of $w(t)$ because the
reasoning is more or less standard. Similar arguments can be found e.g.
either in \cite[Ch.~18, Lemma 24]{M} or in \cite[Th.~3.1]{Y00}.
\cbdu

\medskip

The goal of the next proposition is to prove that the bilinear form
$B(\cdot,\cdot)$ defined in \rf{Bf} is continuous on the space
$\X =\C_w ([0, \infty) , \PM^2)$. This fact is well-known and the
proof
appeared for the first time in \cite{LJS} and \cite{CP2}. Here,
however, we repeat
that reasoning because we want to control better all the constants
which appear
in the estimates below.

\begin{proposition}
\label{prop:L-S}
The bilinear operator $B(\cdot,\cdot)$ is continuous on the space
$\X$ defined in \rf{X}.
Hence, there exists a constant $\eta>0$ such that for every $u,v\in
\X$, it
follows
$$
\|B(u,v)\|_\X\leq \eta \|u\|_\X \|v\|_\X.
$$
\end{proposition}

\proof
We do all the calculations in the Fourier variables. Recall that
the constant $\kappa$
is defined in \rf{kappa}.
Using elementary properties of the Fourier transform we obtain
\begin{eqnarray*}
\left|\widehat \P(\xi) \widehat{(u\otimes v)}(\xi,\tau)\right|
&\leq&
\kappa \int_{\bbfR^3} {dz \over |\xi-z|^2 |z|^2} \;
\|u(\tau)\|_{\PM^2}\|v(\tau)\|_{\PM^2}\\
&=& {\eta \over |\xi|} \|u(\tau)\|_{\PM^2}\|v(\tau)\|_{\PM^2}.
\end{eqnarray*}
In the computations above, we use the equality
$|\xi|^{-2}*|\xi|^{-2}=\pi^{3}|\xi|^{-1}$. A detailed analysis
concerning such convolutions can be found in \cite[Th. 5.9]{LiLo} or
\cite[Ch. V, Sec.1, (8)]{St},
see also \cite[Lem. 2.1]{BCGK}.
Hence,  $\eta=\kappa \pi^3$.

Now, the boundedness of the bilinear form on $\X$ results from the
following
estimates
\begin{eqnarray*}
&& \hspace{-2cm} |\xi|^2 \left|\int_0^t
e^{-(t-\tau)|\xi|^2}\widehat\P (\xi)
\; i\xi  \cdot \widehat{(u\otimes v)} (\xi, \tau)\right|\; d\tau\\
&\leq& \eta |\xi|^2
\int_0^t e^{-(t-\tau)|\xi|^2} \;d\tau \; \|u\|_\X \|v\|_\X \\
&\leq& \eta \|u\|_\X \|v\|_\X.
\end{eqnarray*}
It remains to show the weak continuity of $B(u,v)(t)$
with respect to $t$, but this  follows again from standard
arguments, cf. the remark at the end of the proof of Lemma
\ref{lem-b}.

\cbdu
\medskip

Now, the main theorem of this section results immediately from Lemma
\ref{lem:xyB}
combined with  Lemmata \ref{lem-a}--\ref{lem-b} and Proposition
\ref{prop:L-S}.

\begin{theorem}
Assume that $u_0\in\PM^2$ and $F\in \C_w([0,\infty),\PM)$
satisfy $\|u_0\|_{\PM^2}+\|F\|_{\C_w([0,\infty),\PM)}<\varepsilon$
for some $0<\varepsilon<1/(4\eta)$ where $\eta$ is defined in
Proposition \ref{prop:L-S}.
There exists a global-in-time solution of
\rf{eq}--\rf{ini}  in the space
$
\X= \C_w([0,\infty),\PM^2).
$
This is the unique solution satisfying the condition
$ \|u\|_{\C_w([0,\infty),\PM^2)}\leq 2 \varepsilon.$
Moreover, this solution depends continuously on initial data
and external forces in the sense of Lemma \ref{lem:xyB}. \cbdu
\label{th:glob}
\end{theorem}

Assume, for a moment, that $F\equiv 0$.
Homogeneity properties of the problem \rf{eq}--\rf{div} imply that if
$u$ solves the Cauchy problem, then the rescaled function
$u_\lambda(x,t)=\lambda u(\lambda x,\lambda^2t)$ is also a solution
for each $\lambda>0$. Thus, it is natural to consider solutions which
satisfy the scaling invariance property $u_\lambda\equiv u$ for all
$\lambda>0$, i.e. {\sl forward self-similar } solutions.
By the very definition, they are global-in-time, and one may expect
that they describe the large time behavior of general solutions of
\rf{eq}--\rf{ini}. Indeed, if
$\lim_{\lambda\to\infty}\lambda u(\lambda x,\lambda^2t)=U(x,t)$ in
an appropriate sense, then $tu(xt^{1/2},t)\to U(x,1)$ as $t\to\infty$
(take $t=1$, $\lambda=t^{1/2}$), and $U\equiv U_\lambda$ is scale
invariant. Hence $U$ is a self-similar solution, and
\begin{equation}
U(x,t)=t^{-1/2}U(x/t^{1/2},1)\label{ss}
\end{equation}
is thus determined by a function of $d$ variables $U(y)\equiv
U(y,1)$, $y=x/t^{1/2}$ being the Boltzmann substitution.

If $u_\lambda\equiv u$ for all $\lambda>0$, then from the
self-similar form \rf{ss}, the initial condition \rf{ini}
$\lim_{t\searrow 0}
u(x,t)$ is a distribution homogeneous of degree $-1$
at the origin. Of course, one-point singular solutions defined in
\rf{sing-sol}
are self-similar solutions which are time independent.

Self-similar solutions can be obtained directly from Theorem
\ref{th:glob}
by taking $u_0$ homogeneous of degree $-1$ of small $\PM^{2}$
norm. By the uniqueness property of solutions of the Cauchy problem
constructed in Theorem \ref{th:glob}, they have the form \rf{ss}.

The same reasoning can be applied to the case when external forces
are present. Indeed,
if the initial datum $u_0$ is homogeneous of degree $-1$ and if the
external force $F(x,t)$ satisfies
\begin{equation}
\lambda^3 F(\lambda x, \lambda^2
t)= F(x,t)\quad \mbox{ for all}\quad   \lambda>0
\label{as:F}
\end{equation}
(here, the scaling is understood in
the distributional sense),  the
solution obtained in Theorem \ref{th:glob} is self-similar. Note
that, in particular, we can take
$$
F(x,t)=F(x)= (b_1\delta_0, b_2\delta_0, b_3\delta_0)
$$
(the multiples of the  Dirac delta) for sufficiently small $|b|$. In
other words, the existence of the solutions
introduced by Tian and Xin  and described in the previous section can
be ensured by the fixed point method for
large values of the parameter $c$ (this is possible because
of the particular expression  of the function $b(c)$ in
\rf{fact.b}). We will clarify this  fact in Section 6.

Proceeding in this way we arrive at

\begin{corollary}\label{2}
Suppose that   the initial condition
$u_0\in\PM^{2}$ is homogeneous of degree $-1$ and
$F\in \C_w([0,\infty), \PM)$ satisfies \rf{as:F}. Let $u_0$ and $F$
satisfy, moreover, the assumptions of Theorem \ref{th:glob}.
The corresponding unique solution constructed in
Theorem~\ref{th:glob}
is  self-similar. \cbdu
\end{corollary}

\begin{remark}
The self-similar solutions constructed in such a way can have
singularities for any time. This is the case, for instance, for the
self-similar (stationary)
solutions by
Landau and Tian and Xin. On the other hand, when using the two norms
approach
of Kato as in \cite{C, CaH},
the self-similar solutions that arise from this construction are {\sl
instantaneously} smoothed out for $t>0$
and the only singularity (of the type $\sim 1/\vert x\vert$)
can be found at $t=0$. We will remark on this important point in
Section 7.
\end{remark}

\begin{remark}An
alternative way to prove the existence of self-similar solutions is
to convert \rf{eq}--\rf{ini} into the integral formulation
\rf{duhamel}
and check that the form $B$ reproduces the scale-invariant form
\rf{ss} of $u$. Thus,
the equation \rf{duhamel} can be solved in a subspace of $\x$ formed
by
self-similar functions, as was done in \cite{C}, \cite{M}.
\end{remark}

\medskip

\medskip

\begin{remark}
The existence and the stability results from this section are closely
related to those from the paper by Yamazaki \cite{Y00} where he
studied the Navier--Stokes system in the weak $L^p$-spaces in an
exterior domain $\Omega$. In those considerations, Yamazaki applied
the Kato algorithm in the space $\C_w([0,\infty),
L^{3,\infty}(\Omega))$ without {\sl a priori} assumptions on the
decay of solutions. Our approach involving the $\PM^2$ space is much
more elementary than that from \cite{Y00}. Moreover, we can treat
more singular external forces, and we obtain a kind of asymptotic
stability of solutions (see the next section).
\end{remark}

\begin{remark}
Solutions to the Navier--Stokes system corresponding to  singular
external forces can also be obtained  from very general results by
Kozono and Yamazaki \cite{KY95} where they use the Sobolev-type
spaces
based on homogeneous Morrey spaces.
Their proof
of existence of stationary solutions relies on the inverse function
theorem and
subtle
estimates of the Stokes operator. Next, they investigate properties
of a perturbation of the Stokes operator and they show resolvent
estimates in the Morrey spaces needed in the proof of stability of
stationary solutions. Here, our space $\PM^2$ is much smaller that
those from \cite{KY95}. Our approach, however, besides its
simplicity, does not require separate reasoning for stationary
solutions and unsteady ones. Moreover, we believe that such an
elementary idea will allow to understand better properties of large
solutions (see Section 8).
\end{remark}

\section{Asymptotic behavior of solutions}
\setcounter{equation}{0}
\setcounter{remark}{0}

In our investigations concerning the large time behavior of
solutions to problem \rf{eq}--\rf{ini}
we need the following improvement of Lemma \ref{lem-b}.

\begin{lemma} \label{F:decay}
Assume that $F\in \C_w([0,\infty), \PM)$ satisfies
    $\lim_{t\to\infty} \|F(t)\|_\PM=0$. Then
$$
\lim_{t\to\infty} \left\| \int_0^t S(t-\tau) \P F(\tau)\;d\tau
\right\|_\DPM =0.
$$
\end{lemma}

\proof
It follows from the definition of the norm $\|\cdot\|_\DPM$ that
\begin{eqnarray*}
\left\| \int_0^t S(t-\tau) \P F(\tau)\;d\tau
\right\|_\DPM &\leq & \kappa \sup_{\xi\in\bbfR^3}
\int_0^t |\xi|^2 e^{-(t-\tau) |\xi|^2}\|F(\tau)\|_\PM\;d\tau\\
&\leq& \kappa \left(
    \sup_{\xi\in\bbfR^3} \int_0^{t/2} ... \;d\tau
+
    \sup_{\xi\in\bbfR^3} \int_{t/2}^t ... \;d\tau
\right).
\end{eqnarray*}

Using the substitution $\xi=w\sqrt{t-\tau}$, we first obtain
\begin{eqnarray*}
    \sup_{\xi\in\bbfR^3} \int_0^{t/2} |\xi|^2 e ^{-(t-\tau)|\xi|^2}
\|F(\tau)\|_\PM \;d\tau
&\leq&
\int_0^{t/2}(t-\tau)^{-1} \sup_{w\in\bbfR^3} |w|^2 e ^{-|w|^2}
\|F(\tau)\|_\PM \;d\tau\\
&\leq& C  \int_0^{t/2}(t-\tau)^{-1}
\|F(\tau)\|_\PM \;d\tau\\
&=& C\int_0^{1/2}(1-s)^{-1}
\|F(ts)\|_\PM \;d\tau.
\end{eqnarray*}
Now, the right-hand side of the above inequality tends to 0 as
$t\to\infty$ by the Lebesgue
Dominated Convergence Theorem.

We estimate the  term containing the integral $\int_{t/2}^t ...
\;d\tau$
in the most direct way by
$$
\left(\sup_{\xi\in\bbfR^3} \int_{t/2}^t |\xi|^2 e ^{-(t-\tau)|\xi|^2}
\;d\tau\right) \sup_{t/2\leq \tau\leq t} \|F(\tau)\|_\PM
\leq C  \sup_{t/2\leq \tau\leq t} \|F(\tau)\|_\PM
\to 0
$$
as $t\to\infty$ by the assumption on $F$.
\cbdu

\begin{theorem}\label{asymp}
Let the assumptions of Theorem \ref{th:glob} hold true.
Assume that $u$ and $v$ are two solutions of \rf{eq}--\rf{ini}
constructed in Theorem \ref{th:glob}
corresponding to the initial conditions $u_0,v_0\in\DPM$ and external
forces
$F,G \in \C_w([0,\infty), \PM)$, respectively.
Suppose that
\begin{equation}
\lim_{t\to\infty}\|S(t)(u_0-v_0)\|_\DPM=0
\quad\mbox{and}\quad
\lim_{t\to\infty} \|F(t)-G(t)\|_\PM=0.
\label{e:u0v0}
\end{equation}
Then
\begin{equation}
\lim_{t\to\infty}\|u(\cdot,t)-v(\cdot,
t)\|_\DPM=0\label{e:utvt}
\end{equation}
holds.
\end{theorem}

This result means that if the difference of the solutions of the heat
equation issued from $u_0$, $v_0$ becomes negligible as
$t\to\infty$ (e.g., if the difference of the initial data $u_0-v_0$
is not too singular) and if $F(t)$ and $G(t)$ have the same
large time asymptotics, the
solutions of the nonlinear problem $u(t)$,
$v(t)$ behave similarly for large times. It can be interpreted as a
kind of asymptotic stability result if the choice of $v_0$ is
restricted to the initial data in a~neighborhood of $u_0$ satisfying
additionally \rf{e:u0v0}. It is easy to verify that the first
condition in
\rf{e:u0v0} is satisfied if, e.g.,  $|\xi|^{2}(\widehat u_0(\xi)-\widehat
v_0(\xi))\to 0$ as $\xi\to 0$.
\medskip

\noindent
{\it Proof of Theorem \ref{asymp}.}
First, let us recall that, by Theorem \ref{th:glob}, we have
\begin{equation}
\sup_{t\geq 0} \|u(t)\|_\DPM \leq 2\varepsilon<{1\over 2\eta} \quad
\mbox{and} \quad
    \sup_{t\geq 0} \|v(t)\|_\DPM\leq 2\varepsilon<{1\over 2\eta}.
\label{5.1.u.v}
\end{equation}
We subtract the integral equation \rf{duhamel} for $v$ from the
analogous expression for $u$. Next, computing the norm
$\|\cdot\|_\DPM$ of the resulting equation and repeating the
calculations
from the proof of Proposition \ref{prop:L-S} we obtain the following
inequality
\begin{eqnarray}
&&\hspace{-1cm}\|u(t)-v(t)\|_\DPM \nonumber\\
&\leq& \|S(t)(u_0-v_0)\|_\DPM +
\left\|\int_0^t S(t-\tau)\P
(F(t)-G(t))\;d\tau\right\|_\DPM\label{as-1}\\
&& + \eta \sup_{\xi\in\bbfR^3} \int_0^{\delta t} |\xi|^2
e ^{-(t-\tau)|\xi|^2} (\|u(\tau)\|_{\DPM}+\|v(\tau)\|_{\PM^2})
\|u(\tau) -v(\tau)\|_\DPM
\;d\tau\nonumber\\
&& +  \eta \sup_{\xi\in\bbfR^3} \int_{\delta t}^t |\xi|^2
e ^{-(t-\tau)|\xi|^2} (\|u(\tau)\|_{\DPM}+\|v(\tau)\|_{\PM^2})
\|u(\tau) -v(\tau)\|_\DPM
\;d\tau.\nonumber
\end{eqnarray}
where small constant $\delta>0$ will be chosen later.

In the term on the right-hand side  of \rf{as-1}
containing the integral $\int_0^{\delta t}
...\;d\tau$, we change the variables $\tau=ts$ and we use the
identity
$$
\sup_{\xi\in\bbfR^3} |\xi|^2 e ^{-(1-s)t |\xi|^2}
=((1-s)t)^{-1}\sup_{w\in\bbfR^3}|w|^2e
^{-|w|^2}=((1-s)t)^{-1} e^{-1}
$$
in order to estimate it by
\begin{eqnarray}
&&  \eta \sup_{\xi\in\bbfR^3} \int_0^{\delta} t|\xi|^2 e
^{-(1-s)t|\xi|^2} \|u(ts)-v(ts)\|_\DPM\; ds \nonumber\\
&&\hspace{1.5cm}\times
\left(\sup_{\tau>0} \|u(\tau)\|_{\PM^2}+\sup_{\tau>0}
\|v(\tau)\|_{\PM^2}
\right)\label{as-2}\\
&& \hspace{1cm}
\leq 4\varepsilon \eta e^{-1}  \int_0^{\delta} (1-s)^{-1}\nonumber
\|u(ts)-v(ts)\|_\DPM\; ds.
\end{eqnarray}

We deal with  the term in \rf{as-1}
    containing  $\int_{\delta t}^t ...\;d\tau$ estimating it directly
by
\begin{eqnarray}
&& \eta \left( \sup_{\xi\in\bbfR^3} \int_{\delta t}^t |\xi|^2
e ^{-(t-\tau)|\xi|^2} \;d\tau \right)\left(
    \sup_{\delta t\leq \tau\leq t} \|u(\tau) -v(\tau)\|_\DPM \right)
    4\varepsilon
\label{as-3}\\
&&\hspace{1cm}
= 4\varepsilon \eta \sup_{\delta t\leq \tau\leq t} \|u(\tau)
-v(\tau)\|_\DPM,
\nonumber
\end{eqnarray}
since $\sup_{\xi\in\bbfR^3} \int_{\delta t}^t |\xi|^2
e ^{-(t-\tau)|\xi|^2} \;d\tau = \sup_{\xi\in\bbfR^3} \left(1-e
^{-(\delta t)|\xi|^2}\right)=1$.

Now, we denote
$$
g(t)= \|S(t)(u_0-v_0)\|_\DPM +
\left\|\int_0^t S(t-\tau)\P (F(t)-G(t))\;d\tau\right\|_\DPM,
$$
and it follows from  the assumptions, \rf{e:u0v0} and Lemma
\ref{F:decay} that
\begin{equation}
g\in L^\infty(0, \infty) \quad \mbox{and} \quad
\lim_{t\to\infty} g(t)=0. \label{as-3a}
\end{equation}
Hence, applying \rf{as-2} and \rf{as-3} to \rf{as-1} we obtain
\begin{eqnarray}
\|u(t)-v(t)\|_\DPM &\leq & g(t) +
4\varepsilon \eta e^{-1}  \int_0^{\delta} (1-s)^{-1}
\|u(ts)-v(ts)\|_\DPM\; ds \label{as-4}\\
&& + 4\varepsilon \eta \sup_{\delta t\leq \tau\leq t} \|u(\tau)
-v(\tau)\|_\DPM
\nonumber
\end{eqnarray}
for all $t>0$.

Next, we put
$$
A=\limsup_{t\to\infty} \|u(t)-v(t)\|_\DPM \equiv
\lim_{k\in \bbfN, k\to\infty} \sup_{t\geq k} \|u(t)-v(t)\|_\DPM.
$$
The number $A$ is nonnegative and finite because both
$u,v\in L^\infty ([0,\infty), \DPM)$, and our claim is to show that
$A=0$.
Here, we apply  the Lebesgue Dominated Convergence Theorem to the
obvious inequality
$$
\sup_{t\geq k} \int_0^{\delta} (1-s)^{-1} \|u(ts)-v(ts)\|_\DPM \;ds
\leq \int_0^{\delta} (1-s)^{-1} \sup_{t\geq k} \|u(ts)-v(ts)\|_\DPM
\;ds,
$$
and we obtain
\begin{equation}
\limsup_{t\to\infty}\int_0^{\delta} (1-s)^{-1} \|u(ts)-v(ts)\|_\DPM
\;ds
\leq A\int_0^{\delta} (1-s)^{-1}\;ds =A\log \left({1\over
1-\delta}\right). \label{as-5}
\end{equation}
Moreover, since
$$
\sup_{t\geq k} \sup_{\delta t\leq\tau\leq t} \|u(\tau)
-v(\tau)\|_\DPM
\leq \sup_{\delta k\leq \tau <\infty} \|u(\tau)-v(\tau)\|_\DPM,
$$
we have
\begin{equation}
\limsup_{t\to\infty}  \sup_{\delta t\leq \tau\leq t} \|u(\tau)
-v(\tau)\|_\DPM\leq A.\label{as-6}
\end{equation}

Finally, computing $\limsup_{t\to\infty}$ of the both sides of
inequality \rf{as-4}, and using \rf{as-3a}, \rf{as-5}, and \rf{as-6}
we get
$$
A\leq \left( 4 \varepsilon \eta e^{-1}\log \left({1\over
1-\delta}\right) +4
\varepsilon\eta\right)A.
$$
Consequently,
it follows that
$
A=\limsup_{t\to\infty}\|u(t)-v(t)\|_\DPM=0
$
because
$$
4 \varepsilon \eta \left(e^{-1}\log \left({1\over 1-\delta}\right)
+1\right) <1,
$$
for $\delta>0$ sufficiently small,
by the assumption of Theorem \ref{th:glob} saying that $0<
\varepsilon <1/(4\eta)$.
This completes the proof of Theorem \ref{asymp}.
\cbdu

\medskip

As a direct consequence the proof of Theorem \ref{asymp},
we have also necessary conditions for \rf{e:utvt} to hold. We formulate this
fact
in the following corollary.

\begin{corollary}\label{or:as}
Assume that $u,v \in \C_w([0,\infty), \PM^2)$ are solutions to system
\rf{eq}--\rf{ini}
corresponding to initial conditions $u_0, v_0\in \PM^2$ and external
forces
$F,G \in \C_w([0,\infty), \PM)$, respectively.
Suppose that
\begin{equation}
\lim_{t\to\infty} \|u(t)-v(t)\|_{\PM^2}=0. \label{c0}
\end{equation}
Then
$$
\lim_{t\to\infty} \left\|S(t)(u_0-v_0)+\int_0^t S(t-\tau)\P
(F(\tau)-G(\tau))\;d\tau\right\|_{\PM^2}=0.
$$
\end{corollary}

\proof
As in the beginning of the proof of Theorem \ref{asymp}, we subtract
the integral equation
\rf{duhamel} for $v$ from the same expression for $u$, and we compute
the $\PM^2$-norm
\begin{eqnarray}
&& \hspace{-1cm}\left\|S(t)(u_0-v_0)+\int_0^t S(t-\tau)\P
(F(\tau)-G(\tau))\;d\tau\right\|_{\PM^2}
\nonumber\\
&& \leq\|u(t)-v(t)\|_{\PM^2}\label{c1}\\
&&\hspace{0.5cm}+\eta \sup_{\xi\in\bbfR^3} \int_0^t |\xi|^2
e^{-(t-\tau)|\xi|^2}
(\|u(\tau)\|_{\PM^2}+\|v(\tau)\|_{\PM^2})
\|u(\tau)-v(\tau)\|_{\PM^2}\;d\tau.\nonumber
\end{eqnarray}
The first term on the right-hand side of \rf{c1} tends to zero as
$t\to\infty$
by  \rf{c0}. To show the decay of  the second one, it suffices to
repeat calculations
from \rf{as-1}, \rf{as-2}, \rf{as-3}, and \rf{as-5}. Here, however,
one should remember that now it is assumed that
$A=0$ and $\sup_{t>0} \|u(t)\|_{\PM^2}<\infty$ and
   $\sup_{t>0} \|v(t)\|_{\PM^2}<\infty$.
   \cbdu

\begin{remark}
The Lyapunov stability of solutions (not necessarily stationary ones)
follows immediately from the construction {\sl via} the Banach fixed
point theorem (cf. Lemma \ref{lem:xyB}). This phenomenon was already
observed and used several times, see e.g. the papers by H. Kozono and
M. Yamazaki \cite[Th. 2]{KY95}, \cite[Th. 1]{KY-Zeit}, and by M.
Yamazaki
\cite[Th. 1.3]{Y00}. Theorem \ref{asymp} extends those results by
giving sufficient conditions on the {\sl asymptotic} stability of
solutions. In particular, Yamazaki \cite[Remark 4.1]{Y00} emphasized
that the trivial solution 0 is stable but not asymptotically stable
in the space $L^{3,\infty}(\bbfR^3)$ (in contrast to the Lebesgue
space
$L^3(\bbfR^3)$), because there exist self-similar solutions with
constant $L^{3,\infty}$-norm.
Theorem \ref{asymp} and Corollary \ref{or:as}  explain this
phenomenon in the case of the space $\PM^2$. Indeed, given $u_0\in\PM^2$ such
that $\|u_0\|_{\PM^2}<\varepsilon$ and $F\equiv 0$, the corresponding
solution converges in $\PM^2$ to zero as $t\to\infty$ if and only if
$\lim_{t\to\infty}\|S(t)u_0\|_{\PM^2}=0$.

Note here, that  if $U(x,t)=t^{-1/2} U(x/t^{1/2})$ is a
self-similar solution to system \rf{eq}--\rf{ini}, its $\PM^2$-norm
is
constant in time by the scaling relation \rf{PMa-scal}. As it is
well-known, $U(x,t)$ corresponds to the initial condition $U_0(x)$
which is homogeneous of degree $-1$, so $S(t)U_0(x)= t^{-1/2}
S(1)U_0(xt^{-1/2})$. Consequently, by the scaling property of the
norm, we have
$
\|S(t)U_0\|_{\PM^2} =\|S(1) U_0\|_{\PM^2},
$
cf. Corollary \ref{or:as} with $F=G\equiv 0$.

\end{remark}

\begin{remark}
In the
setting of the $L^p$-spaces and the homogeneous Besov spaces, the
study of the asymptotic stability of
self-similar solutions to the Navier--Stokes system begun with the
paper \cite{P} of F. Planchon (see
also the presentation of Planchon's results in \cite[Ch.~23.3]{LR}).
As illustrated in the book by Y. Giga and M.-H. Giga \cite{GG}
those ideas are quite universal and were used for other partial
differential equations (e.g. the porous medium, the nonlinear
Schr\"odinger and
the KdV equations); they were
applied for instance to study
asymptotic properties of solutions to a large class of nonlinear
parabolic equations \cite{K} as well as of solutions with zero mass
to viscous conservation laws \cite{KS}. In this section, we extend
them on solutions which not necessarily decay to 0 as
$t\to\infty$.
\end{remark}

\section{Stationary solutions}
\setcounter{equation}{0}
\setcounter{remark}{0}

Our approach, described in previous sections,  to study
global-in-time solutions to the problem \rf{eq}--\rf{ini}, as well as
their
large time behavior, can be also applied to stationary solutions.
Below, we briefly describe some consequences of Theorems
\ref{th:glob} and \ref{asymp}. The following proposition contains two
equivalent
integral equations satisfied by stationary solutions.

\begin{proposition}\label{prop:stat}
Assume that $u=u(x)\in \PM^2$ and $F\in \PM$. The following two facts
are equivalent
\begin{itemize}
\item[1)] $u=u(x)$ is a stationary mild solution of system
\rf{eq}--\rf{div}
in the sense of Definition \ref{def1}. Hence,  $u$ is the
solution of the integral equation
\begin{equation}
u= S(t) u - \int_0^t  S(t-\tau) \P \nabla \cdot (u\otimes
u)\;d\tau+\int_0^t
S(\tau) \P F \;d\tau \label{duh1}
\end{equation}
for every $t>0$;

\item[2)] $u$ satisfies the integral equation
\begin{equation}
u=-\int_0^\infty S(\tau ) \P\nabla (u\otimes u)\;d\tau +\int_0^\infty
S(\tau) \P F \;d\tau, \label{duh2}
\end{equation}
where the integrals above should be understood in the Fourier
variables for almost every $\xi$.

\end{itemize}
\end{proposition}

\proof
By Definition \ref{def1}, the integral equation \rf{duh1} can be
rewritten
as
\begin{eqnarray}
\widehat u(\xi) &=& e^{-t|\xi|^2} \widehat u(\xi) -\int_0^t
e^{-(t-\tau)|\xi|^2} \,d\tau \; \widehat\P (\xi) i\xi \cdot
\widehat{(u\otimes
u)}(\xi) \nonumber\\
&&+\int_0^t e^{-(t-\tau)|\xi|^2} \;d\tau \; \widehat\P (\xi)
\widehat F(\xi)\label{duh3}\\
&=& e^{-t|\xi|^2} \widehat u(\xi) -{1-e^{-t|\xi|^2}\over |\xi|^2}
\widehat\P (\xi) i\xi\cdot \widehat{(u\otimes u)}(\xi)
+{1-e^{-t|\xi|^2}\over |\xi|^2} \widehat\P (\xi)
\widehat F(\xi).\nonumber
\end{eqnarray}
for every $t>0$. Passing to the limit as $t\to\infty$ in \rf{duh3}
and using the
identity
$$
{1\over |\xi|^2} = \int_0^\infty e^{-\tau |\xi|^2}\;d\tau \quad
\mbox{for $\xi\neq 0$},
$$
we obtain equation \rf{duh2} in the Fourier variables.

Now, assume that $u$ solves \rf{duh2}. Repeating the arguments above
in the reverse order, we obtain that $u$ is the solution of the
equation
\begin{equation}
\widehat u(\xi)=
-{1\over |\xi|^2}
\widehat\P (\xi) i\xi\cdot  \widehat{(u\otimes u)}(\xi)
+{1\over |\xi|^2} \widehat\P (\xi)
\widehat F(\xi).\label{duh4}
\end{equation}
If we subtract from this equality the same expression multiplied by
$e^{-t|\xi|^2}$ we get \rf{duh3} which obviously is equivalent to
\rf{duh1}.
\cbdu

\begin{theorem}\label{th:stat}
Assume that $F\in \PM$ satisfies $\|F\|_\PM<\varepsilon<1/(4\eta)$. There exists
a~stationary solution $u_\infty$ to the Navier--Stokes system in the
space
$\PM^2$ with $F$ as the external force. This is the unique solution
satisfying the condition
$\|u\|_{\PM^2}\leq 2\varepsilon$.
\end{theorem}

\proof
This theorem results immediately from Lemma \ref{lem:xyB} applied to
the integral equation \rf{duh2} (or its equivalent version
\rf{duh4}). The bilinear form
$$
B(u,v)= \int_0^\infty S(\tau)
\P\nabla\cdot u\otimes v\;d\tau
$$
is bounded on the space $\PM^2$ and
the proof of this property of $B(\cdot,\cdot)$ is completely
analogous to the
one of Proposition \ref{prop:L-S}. Let us also skip an easy proof
that
$y= \int_0^\infty S(\tau) \P F \; d\tau$ satisfies $\|y\|_{\PM^2}
=\|F\|_\PM$.
\cbdu

\medskip

Now, the application of Theorem \ref{asymp} gives the following
result on the asymptotic stability of stationary solutions.

\begin{corollary}
Assume that $u_\infty$ is the stationary solution constructed in
Theorem \ref{th:stat} corresponding to the external force $F$.
Suppose that $v_0\in \PM^2$ and $G\in \C_w([0,\infty), \PM)$ satisfy
$\|v_0\|_{\PM^2}+ \|G\|_{\C_w([0,\infty), \PM)}\leq
\varepsilon<1/(4\eta)$ and, moreover,
$$
\lim_{t\to\infty}\|S(t)(v_0-u_\infty)\|_{\PM^2}=0, \quad
\lim_{t\to\infty}\|G(t)-F\|_{\PM}=0.
$$
Then, the solution $v=v(x,t)$ of system \rf{eq}--\rf{ini}
corresponding
to $v_0$ and $G$  converges toward the stationary solution
$u_\infty$ in the following sense
$$
\lim_{t\to\infty} \|v(t)-u_\infty\|_{\PM^2}=0.
$$
\end{corollary}

\proof
Here, it suffices only to note that stationary solutions belong to
the space $\C_w([0,\infty), \PM^2)$ (treated as constant functions on
$[0, \infty)$ with values in $\PM^2$) and satisfy the integral
equation \rf{duhamel} (see Proposition \ref{prop:stat}). So, Theorem
\ref{asymp} is applicable in this
case. \cbdu

\begin{remark}
Results from this section can be extended to solutions which exist
for all $t\in\bbfR$ (and not only for $t\geq 0$) as  was done by
M. Yamazaki \cite{Y00}. In this case the corresponding integral
equation
(the counterpart of \rf{duh1} and \rf{duh2}) has the form
$$
u(t)=-\int_0^\infty S(\tau ) \P\nabla (u\otimes u)(t-\tau)\;d\tau
+\int_0^\infty
S(\tau) \P F (t-\tau)\;d\tau,
$$
and, like in \cite{Y00}, by the application of Theorem \ref{th:glob},
one obtains solutions which  are, for example, time periodic or
almost periodic with respect to $t\in\bbfR$. In the same manner,
Theorem \ref{asymp} allows us to describe solutions which converge in
$\PM^2$ as $t\to\infty$ toward given time periodic (or
almost periodic) solution.
\end{remark}


\section{Smooth solutions}
\setcounter{equation}{0}
\setcounter{remark}{0}

Solutions of problem \rf{eq}--\rf{ini} constructed in the space
$\X=\C_w([0, \infty), \PM^2)$ are, in fact, smooth (for sufficiently
regular external forces), and they agree with mild solutions obtained
by T. Kato \cite{Ka} and in \cite{C} for $F\equiv 0$, and,
more generally,
with solutions obtained in \cite{CPfo} when $F\not =0$.

The goal of this section is to clarify this remark.
First, let us recall that, in \cite{C}, solutions of
\rf{eq}--\rf{ini}
were constructed for sufficiently small initial conditions from the
homogeneous Besov space
${\dot B}^{-1+3/p,\infty}_p(\bbfR^3)$ with $3<p<\infty$. The usual
way of defining a norm in this space is based on the dyadic
decomposition of tempered distributions. Here, however as in 
\cite{C,K}
we prefer the
equivalent norm whose definition involves the heat semigroup
$$
\|v\|_{{\dot B}^{-\alpha,\infty}_p(\bbfR^3)}\equiv \sup_{t>0}
t^{\alpha/2}
\|S(t)v\|_{L^p(\bbfR^3)}.
$$
Connections between $\PM^2$ and homogeneous Besov spaces are
described in the following lemma.

\begin{lemma}
For every $p\in (3, \infty]$ the following imbeddings
$
\DPM \subset {\dot B}^{-1+3/p,\infty}_p(\bbfR^3)
$
hold true and are continuous.
Hence, there exists a constant $C=C(p)$  such that
$$
\sup_{t>0}t^{(1-3/p)/2}\|S(t)u_0\|_{L^p(\bbfR^3)} \leq C
\|u_0\|_\DPM
$$
for all $t>0$ and $u_0\in \PM^2$.
\end{lemma}

\proof
Here, our tool is the Hausdorff--Young inequality. For $1/p+1/q=1$ we
obtain
\begin{eqnarray*}
\|S(t)u_0\|_{L^p(\bbfR^3)}^q &\leq & C \int_{\bbfR^3}
\left|e^{-t|\xi|^2} \widehat
u_0(\xi) \right|^q\,d\xi\\
&\leq& C \sup_{\xi\in\bbfR^3} |\xi|^2 |\widehat u_0(\xi)|
\int_{\bbfR^3} \left| {e^{-qt|\xi|^2}\over |\xi|^{2q}}\right|\,d\xi\\
&=& C\|u_0\|_\DPM t^{-3/2+q}
\int_{\bbfR^3} \left| {e^{-q|w|^2}\over |w|^{2q}}\right|\,dw.
\end{eqnarray*}
In the calculations above, we assume that $2q<3$ which is equivalent
to $p>3$. Since, $1/q=1-1/p$ and $(3/2)(1-1/p)-1=(1/2)(1-3/p)$, we
obtain
$$
t^{(1/2)(1-3/p)}\|S(t)u_0\|_{L^p(\bbfR^3)}\leq C\|u_0\|_\DPM.
$$

Note that this proof requires an obvious modification
for $p=\infty$ and $q=1$. One can also  recall here  the embedding of
any ``critical space'' into the
Besov space ${\dot B}^{-1,\infty}_\infty(\bbfR^3)$, see \cite{M,
CaH}.
\cbdu

\medskip

Now, given $u_0\in\PM^2$ with sufficiently small $\PM^2$-norm, we may
apply the theory described in \cite{C}  to get the solution
$\widetilde u =\widetilde u(x,t)$ which is unique in the space
$$
\C_w([0, \infty), {\dot B}^{-1+3/p,\infty}_p(\bbfR^3)) \cap
\{v\;:\; t^{(3/p-1)/2}\|v(t)\|_{L^p(\bbfR^3)}<\infty\}
$$
corresponding to $u_0$ as the initial condition and the zero external
force. Moreover, this solution is smooth for all $t>0$.
On the other hand, our Theorem \ref{th:glob} gives a solution
$u=u(x,t)$ in $\C_w([0,\infty), \PM^2)$.

Both constructions lead, in fact, to the same solution, and we show
this by analyzing
    the parabolic
regularization effect in problem \rf{eq}--\rf{ini}
in the scale of spaces $\PMa$. We begin by a definition.

\begin{definition}\label{def2}
Let  $2\leq a< 3$. We define the Banach space
\begin{eqnarray}
\ya&\equiv& \C_w([0, \infty), \PM^2)\label{Ya}\\
&&\cap\;\; \{v:(0,\infty)\to
\PM^a:
\I v\I_a\equiv\sup_{t>0}t^{a/2-1}\|v(t)\|_{\PM^a}<\infty\}.\nonumber
\end{eqnarray}
The space $\ya$ is normed by the quantity
$\|v\|_{\ya}=\I v\I_{2}+\I v\I_a$.
Of course, ${\cal Y}^{2}\equiv \X$ with this definition.
\end{definition}

\begin{remark}
The norm $\I\cdot \I_a$ is invariant under the rescaling
$u_\lambda(x,t)=\lambda
u(\lambda x, \lambda^2 t)$ for every $\lambda >0$. This can be easily
checked using the scaling property of the norm $\|\cdot\|_{\PMa}$,
see \rf{PMa-scal}.
\end{remark}
\medskip

First we show an improvement of Proposition \ref{prop:L-S}.

\begin{proposition}\label{prop:7.1}
Let $2\leq a<3$. There exists a constant $\eta_a>0$ such that
for every $u\in \C_w ([0,\infty), \PM^2)$ and
$v\in \{v(t)\in\PM^a\;:\; \I v\I_a<\infty\}$ we have
$$
\I B(u,v)\I_a\leq \eta_a \I u\I_2 \I v\I_a.
$$
\end{proposition}

\proof
First note that as in the proof of Proposition \ref{prop:L-S} we have
\begin{eqnarray*}
|\widehat{(u\otimes v)}(\xi,t)|&\leq& \int_{\bbfR^3} {1\over
|\xi-z|^2|z|^a}\;dz \;\|u(t)\|_{\PM^2}\|v(t)\|_\PMa\\
&=& C|\xi|^{1-a} \|u(t)\|_{\PM^2}\|v(t)\|_\PMa.
\end{eqnarray*}
Thus, for every $\xi\neq 0$ we obtain
$$
|\xi|^a \left|\int_0^t e^{-(t-\tau)|\xi|^2} \widehat\P (\xi)
    i\xi\cdot \widehat{(u\otimes v)}(\xi,\tau)\;d\tau\right|
\leq C \int_0^t |\xi|^2 e^{-(t-\tau)|\xi|^2} \tau^{1-a/2}\;d\tau
\;\I u\I_2\I v \I_a.
$$

The proof will be completed by showing that for every $a\geq 2$ the
quantity
\begin{equation}
t^{a/2-1}\int_0^t |\xi|^2 e^{-(t-\tau)|\xi|^2} \tau^{1-a/2}\;d\tau
\label{t-a-xi}
\end{equation}
is bounded by a constant independent of $\xi$ and $t$.
Here, we decompose the integral with respect to $\tau$ into two parts
$
\int_0^t...\;d\tau = \int_0^{t/2}...\;d\tau +\int_{t/2}^t...\;d\tau,
$
and we deal with the each term separately.

In case of  the integral over $[0,t/2]$, we estimate the above
quantity by
$$
t^{a/2-1} |\xi|^2 e^{-(t/2)|\xi|^2}\int_0^{t/2} \tau^{1-a/2}\;ds
= C(t/2)|\xi|^2e^{-(t/2)|\xi|^2} \leq C
$$
where $C$ is independent of $\xi$ and $t$. For the interval
$[t/2,t]$, the quantity is bounded by
$$
t^{a/2-1} (t/2)^{1-a/2} \int_{t/2}^t |\xi|^2
e^{-(t-\tau)|\xi|^2}\;d\tau
=(1/2)^{1-a/2} (1-e^{-(t/2)|\xi|^2}) \leq (1/2)^{1-a/2}.
$$
\cbdu
\medskip

Next, we show that that the heat semigroup regularizes distributions
from $\PM^2$.

\begin{lemma}\label{lem:S(t)}
For every $u_0\in \PM^2$ and $t>0$, it follows that $S(t)u_0\in
\PM^a$ with $a\geq 2$. Moreover, there exists $C$ depending on 
the exponent $a$
only such that
$$
\sup_{t>0} \left(t^{a/2-1} \|S(t)u_0\|_{\PMa}\right) \leq C
\|u_0\|_{\PM^2}.
$$
\end{lemma}

\proof
Simple estimates (cf. Lemma \ref{lem-a}) give
\begin{eqnarray*}
\sup_{t>0} \left(t^{a/2-1} \|S(t)u_0\|_{\PMa}\right) &\leq &
   \|u_0\|_{\PM^2} \sup_{\xi\in\bbfR^3}
\left(t^{a/2-1}|\xi|^{a-2} e^{-t|\xi|^2}\right)\\
&=& C \|u_0\|_{\PM^2}
\end{eqnarray*}
where $C= \sup_{w\in\bbfR^3}
\left(|w|^{a-2} e^{-|w|^2}\right)$.
\cbdu

\medskip
Let us also explain how to handle more regular external forces in the
scale of the spaces $\PMa$.

\begin{lemma}\label{7.3}
Let $2\leq a<3$. Assume that $F(t)\in\PM^{a-2}$ for all $t>0$ and
\begin{equation}
\sup_{t>0} t^{a/2-1}\|F(t)\|_{\PM^{a-2}}<\infty.\label{F:PMa}
\end{equation}
There exists a
constant $C$ such that for $w(t)=\int_0^t S(t-\tau)\P F(\tau)\;d\tau$
it follows that
$$
\I w\I_a\leq C \sup_{t>0} t^{a/2-1}\|F(t)\|_{\PM^{a-2}}.
$$
\end{lemma}

\proof
As in the proof of Lemma \ref{lem:S(t)}, we obtain
\begin{eqnarray*}
\|w(t)\|_{\PM^a} &\leq & \esssup |\xi|^a \int_0^t
\left| e^{-(t-\tau)|\xi|^2} \widehat \P(\xi) \widehat F
(\xi,\tau)\right|\;d\tau\\
&\leq& \kappa \int_0^t |\xi|^2 e^{-(t-\tau)|\xi|^2}
\tau^{1-a/2}\;d\tau \; \sup_{t>0} t^{a/2-1}\|F(t)\|_{\PM^{a-2}}.
\end{eqnarray*}
  From now on, it suffices to repeat the reasoning which leads to the
estimates of the quantity in \rf{t-a-xi}.
\cbdu

\begin{theorem}\label{th:reg}
Let $a\in [2,3)$.  There exists
$\varepsilon>0$ such that for
every $u_0\in \PM^2$ and $F\in \C_w([0,\infty), \PM)$ satisfying
\rf{F:PMa} with
$$
\|u_0\|_{\PM^2}+\|F\|_{\C_w([0,\infty), \PM)}+
\sup_{t>0} t^{a/2-1}\|F(t)\|_{\PM^{a-2}}
<\varepsilon,
$$
the solution
constructed in Theorem \ref{th:glob} satisfies
$
\I u\I_a\leq 2\varepsilon.
$
\end{theorem}

\proof
It suffices to repeat the reasoning leading to  Theorem
\ref{th:glob} in the space
$$
{\cal Y}^a= \C_w([0,\infty), \PM^2)\cap \{u\;:\; \sup_{t>0}
t^{a/2-1}\|u(t)\|_{\PMa}<\infty\}
$$
involving Lemma \ref{lem:xyB}. Here, the required estimate of the
bilinear form $B(\cdot,\cdot)$ is proved
in Propositions \ref{prop:L-S} and \ref{prop:7.1}.
Moreover, Lemmata \ref{lem:S(t)} and \ref{7.3} guarantee that
$y=S(t)u_0 +\int_0^t S(t-\tau)\P F(\tau)\;d\tau$
belongs to ${\cal Y}^a$.
\cbdu

\medskip

Let us formulate an interpolation inequality involving $L^q$ and $\PMa$ norms.

\begin{lemma}\label{lem:interpol}
Fix $a\in (2,3)$. For every $q\in \left(3,{3\over 3-a}\right)$ there
exists a
constant $C=C(a,q)$  such that
\begin{equation}
\|v\|_{L^q(\bbfR^3)} \leq C\|v\|_{\PM^2}^\beta
\|v\|_{\PMa}^{1-\beta}
\label{interpol}
\end{equation}
for all $v\in \PM^2\cap\PMa$, where $\beta= {1\over
a-2}\left(1-{3\over q}\right)$.
\end{lemma}

\proof
Assume that $v$ is smooth and rapidly decreasing. Using the
Hausdorff--Young inequality (with $1/p+1/q=1$ and $p\in [1,2)$) and
the definition of the $\PMa$-norm we obtain
\begin{eqnarray}
\|v\|_q^p&\leq& C\|\widehat v \|_p^p \leq C\|v\|_2^p \int_{|\xi|\leq
R}
{1\over |\xi|^{2p}}\;d\xi +C\|v\|_a^p \int_{|\xi|> R}
{1\over |\xi|^{ap}}\;d\xi\nonumber\\
&\leq& C\|v\|_2^pR^{3-2p}+C\|v\|_a^pR^{3-ap}\label{interpol2}
\end{eqnarray}
for all $R>0$ and $C$ independent of $v$ and $R$. In these
calculations, we require $2p<3$ which is equivalent to $q>3$.
Moreover, we have to assume that $ap>3$ which leads to the inequality
$q<3/(3-a)$. Now, we optimize inequality \rf{interpol2} with respect
to $R$ to get \rf{interpol}.
\cbdu

\begin{corollary}\label{cor:7.1}
Under the assumptions of Theorem \ref{th:reg} the constructed
solution satisfies
$$
\|u(\cdot,t)\|_{L^q(\bbfR^3)}\leq Ct^{-(1-3/q)/2}
$$
for each $3<q,3/(3-a)$, all $t>0$, and $C$ independent of $t$.
\end{corollary}

\proof
It follows from Theorem \ref{th:reg} that the solution $u$ satisfies
$\|u(\cdot,t)\|_{\PMa}\leq C t^{1-a/2}$ for every $a\in [2,3)$.
Hence, to complete the proof of this corollary, it suffices to apply
Lemma \ref{lem:interpol}. \cbdu

\bigskip

Let us finally prove that the difference of two (singular) solutions
corresponding to the same external
force is more regular than each term separately. This
fact is in a perfect
agreement with the regularity result for the bilinear term obtained
in \cite{CP2}.

\begin{theorem}
Assume that $u,v\in\X$ are solutions to \rf{eq}--\rf{ini}
constructed in Theorem \ref{th:glob} corresponding to initial
conditions $u_0,v_0\in \PM^2$ and the same external force
$F\in \C_w([0, \infty), \PM)$.
For every $2\leq a<3$ there exists $\varepsilon>0$ such that
for $\|u_0-v_0\|_{\PM^2}<\varepsilon$ we have
$$
\I u-v\I_a \equiv \sup_{t>0} t^{a/2-1} \|u(t)-v(t)\|_{\PM^a}<\infty.
$$
Moreover, $\sup_{t>0} t^{(1-3/q)/2}
\|u(t)-v(t)\|_{L^q(\bbfR^3)}<\infty$ for every $3<q<3/(3-a)$.
\end{theorem}

\proof
Here, the reasoning is similar to that presented above, hence we
shall be brief in details. First, we subtract integral equations
\rf{duhamel} for $u$ and $v$ to obtain
$$
u(t)-v(t)=S(t)(u_0-v_0) +B(u,u-v)(t)+B(u-v,v)(t).
$$
We denote $z(t)=u(t)-v(t)$ and $z_0=u_0-v_0$, and we find the
solution of the equation
$
z=S(\cdot)z_0+B(u,z)+B(z,v)
$
{\sl via} the Banach fixed point theorem in the space $\ya$ defined
in \rf{Ya}. Here, Lemma \ref{lem:S(t)} guarantees that
$S(\cdot)z_0\in \ya$ for every $2\leq a<3$. Moreover, Propositions
\ref{prop:L-S} and \ref{prop:7.1} allow us to show the contractivity
of the mapping $z\mapsto S(\cdot)z_0+B(u,z)+B(z,v)$ for sufficiently
small $\varepsilon>0$ because, by Theorem \ref{th:glob}, $u$ and $v$
satisfy \rf{5.1.u.v}. The second part of this theorem is deduced
immediately from Lemma \ref{lem:interpol}.
\cbdu

\medskip

\begin{remark}
Given $u_0\in \PM^2$ with sufficiently small norm and $F\equiv
0$, Theorem \ref{th:glob} guarantees the  existence of a  unique
small
solution $u\in \C_w([0,\infty), \PM^2)$. Next, our analysis in
Corollary \ref{cor:7.1} allows us to show that $u(t)\in L^q(\bbfR^3)$
for  $q>3$ and all $t>0$. Hence, standard regularity theorems
imply that $u(x,t)$ is a smooth function and satisfies the
Navier--Stokes system in the classical sense.

Even if it is not written explicitly, the same conclusion can be
deduced from Yamazaki's results
\cite[Th.~1.3]{Y00}, where he
showed that his solution belonging initially to $\C_w([0,\infty),
L^{3,\infty} (\bbfR^3))$ falls, in fact, into $L^{p,\infty}(\bbfR^3)$
for every $3<p<\infty$. Now, applying the Marcinkiewicz interpolation
theorem for the identity mapping, one obtains immediately that
$
\bigcap_{3<p<\infty} L^{p,\infty}(\bbfR^3) \subset L^q(\bbfR^3)
$
for every $3<q<\infty$.
\end{remark}

\medskip
We conclude this section by
stressing again that the two norms approach by Kato imposes {\sl
a~priori} a~regularization effect on solutions we look for. In other
words, they are considered as fluctuations around the solution of the
heat
equation $S(t)u_0$. The solutions appear to be unique locally
in the space of more regular functions. The approach with the only
one norm in Theorem \ref{th:glob} gives the local uniqueness in the
larger
space which, in our case, may contain genuinely singular solutions
(like those in \rf{sing-sol}) which are not smoothed out by the
action of the nonlinear
semigroup associated with \rf{eq}--\rf{ini}.


\section{Loss of smoothness for large solutions}
\setcounter{equation}{0}
\setcounter{remark}{0}

As far as blow-up for Navier--Stokes several possibilities can be
conjectured. One may imagine that blow-up of
initially regular solutions never happens, or  it becomes more likely
as the initial norm increases, or that
there is blow-up, but only on a very thin set, of
measure zero.

As we have seen in the previous sections, when using a fixed point
approach,
existence and uniqueness of global solutions are guaranteed
only under restrictive   assumptions on the initial data and external
forces, that
are required to be small in some sense, i.e.  in some functional
space. In \cite{C} we pointed out that  fast oscillations are
sufficient to  make the fixed point scheme works, even if the
norm in the corresponding function space of the initial data is
arbitrarily large (in fact, a different auxiliary norm turns out
to be small).  Here we want to suggest how some particular
data, arbitrarily large (not oscillating) could give rise to
irregular solutions. It is extremely
unpleasant  that we do not know in general whether for arbitrary
large
data the corresponding solution is  regular or singular. More precisely:

\begin{remark} Let us consider the Navier--Stokes equations \rf{ini}
with external force $F\equiv 0$. Then, if one defines the functions
$u_\varepsilon(x,0)=\varepsilon
u(x)$, where $u(x)$ is the (divergence free, homogeneous of degree
$-1$) function given by
\rf{sing-sol} as the initial data, then  for small $\varepsilon$ the
system has a global regular (self-similar) solution
which is even more regular than {\sl a priori} expected (Section 7)
and for $\varepsilon=1$ (and possibly for other large values of
$\varepsilon$)
the system has a~singular ``solution'' for any time.
\end{remark}

\medskip

Unfortunately, this loss of smoothness for large data does not hold
in the
``distributional'' sense,
but as explained in Section 2, only ``pointwise'' for  every $x\in
\bbfR^3\setminus
\{(0,0,0)\}$.
However, for a model equation of gravitating particles
this loss of smoothness for large data holds in the distributional
sense and will be dealt with in a forthcoming
paper \cite{BCGK}.


\bigskip
\noindent
{\bf Acknowledgements.}~The preparation of this paper was partially
supported by the KBN
grant 50/P03/2000/18 and the POLONIUM project \'EGIDE--KBN
2002--2003.
The authors acknowledge gratefully fruitful discussions and  helpful
remarks of Piotr Biler and
Andrzej Krzywicki, in particular concerning calculations in Section
2.

\end{document}